\documentclass[portuguese,a4paper,twoside,10pt]{article}
\usepackage[english]{babel}
\usepackage[T1]{fontenc}
\usepackage[latin1]{inputenc}
\usepackage{amssymb,amsmath,epsfig,amsfonts,euscript}
\usepackage{dsfont}
\usepackage[usenames,dvipsnames]{color}

\newtheorem{teo}{Theorem}[section]
\newtheorem{pro}[teo]{Proposition}

\newcounter{note}[section]
\newcounter{example}[section]

\newcommand{\dem}{{\bf Proof. }}
\newcommand{\fdem}{ $\square$}

\def\rain{\to \infty}
\def\conv{\hspace{0.25em}{\atop{\scriptstyle
n\rain}}\hspace{-2.5em} \longrightarrow \hspace{0.2cm}}

\def\N{{\rm I\kern-.20em N}}
\def\R{{\rm I\kern-.20em R}}
\def\indi{{1\kern-.20em\rm I}}

\newcommand{\titulo}[1]{\begin{center}\mbox{} \\ \noindent \textit{\textbf{\Large #1}}\\\vspace{0.5cm}\end{center}}

\newcommand{\autor}[1]{\noindent \textbf{#1}}

\newcommand{\afil}[1]{{\small \noindent \textit{#1}}}

\renewcommand{\abstract}[1]{{\small \noindent \textbf{Abstract:} #1\\}}

\newcommand{\keywords}[1]{{\small \noindent \textbf{Keywords:} #1\\}}

\newcommand{\biblio}[1]{{\small \begin{thebibliography}{99} #1 \end{thebibliography}}}

\begin{document}
\titulo{Multivariate maxima of moving multivariate maxima}

\autor{Helena Ferreira \footnote{E-mail:helenaf@ubi.pt }} \afil{Department of Mathematics, University
of Beira Interior, Covilh\~{a},
Portugal}\\

\abstract{We define a class of multivariate maxima of moving multivariate maxima, generalising the M4 processes. For these stationary multivariate time series we characterise the joint distribution of extremes and compute the multivariate extremal index. We derive the bivariate upper tail dependence coefficients and the extremal coefficient of the new limiting multivariate extreme value distributions.}

\section{Introduction}


The M$4$ processes considered by Smith and Weissman (1996)
constitute quite a large and flexible class of multivariate time series models exhibiting clustering of extremes. In the following generalisation we will introduce a copula in the array of innovations, rendering the new model more flexible than the previously mentioned M$4$ processes, to which the new model reduces if the copula is  the Fr\'echet-Hoeffding upper bound. The new family (M$5$) comprises multivariate time series models exhibiting clustering of extremes across variables and over time. 
We derive the limiting distribution of the vector of componentwise maxima and the corresponding limit for the sequence of independent and identically distributed variables.
Because of temporal dependence, these two limiting distributions are different, and the difference is quantified by the multivariate extremal index, which will be calculated too. \\

Let $\{{\bf Z}_{l,n}=(Z_{l,n,1},\ldots,Z_{l,n,d})\}_{l\geq 1
,-\infty<n<\infty}$ be an array of independent random vectors with standard Fr\'echet margins and common copula $C_{\bf Z}$. A multivariate maxima of moving multivariate maxima (henceforth M$5$) process is defined by 
\begin{eqnarray}
Y_{n,j}=\max_{l\geq
1}\max_{-\infty<k<+\infty}\alpha_{l,k,j}Z_{l,n-k,j},\,\,j=1,\ldots,d,\,\,n\geq
1,\label{1}
\end{eqnarray}\vspace{0.5cm}
\noindent where $\{\alpha_{l,k,j}, l\geq 1, -\infty<k<\infty,1\leq
j\leq d\}$ are nonnegative constants satisfying
\begin{eqnarray*}
\sum_{l=1}^\infty
\sum_{k=-\infty}^\infty \alpha_{l,k,j}=1,\quad {\rm for}\
j=1,\ldots,d.
\end{eqnarray*}
When $Z_{l,n,j}=Z_{l,n}$, $j=1,...,d$, the M$5$ process is the M$4$ process considered by Smith and Weissman (1996).

The common distribution $F_{\bf Y}$ of ${\bf Y}_{n}=(Y_{n,1},\ldots,Y_{n,d})$ satisfies 
{\small{\begin{eqnarray*}
F_{{\bf Y}}(y_1,...,y_d)=
\prod_{l=1}^{\infty}\prod_{k=-\infty}^{\infty}F_{\bf Z}\left(\frac{y_1}{\alpha_{l,k,1}},...,\frac{y_d}{\alpha_{l,k,d}}\right),\,\,\, 
y_j>0,\, j=1,...,d,
\end{eqnarray*}}}\vspace{0.5cm}
and  the relation for the corresponding copulas is 
{\small{\begin{eqnarray*}
C_{{\bf Y}}(u_1,...,u_d)=
\prod_{l=1}^{\infty}\prod_{k=-\infty}^{\infty}C_{\bf Z}\left(u_1^{\alpha_{l,k,1}},...,u_d^{\alpha_{l,k,d}}\right),\,\,\, 
u_j \in [0,1],\, j=1,...,d.
\end{eqnarray*}}}\vspace{0.5cm}

This paper is concerned with the extreme value properties of the class of stationary processes defined in (\ref{1}). The M5 processes contribute to the modelling  of variables moving together and generalise the important class of moving maxima discussed by several authors
(Davis and Resnick, 1989; Deheuvels, 1983; Hall {\it et al.}, 2002; Smith and Weissman, 1996; Zhang and Smith, 2004; among others).
 This is the motivation of this paper which is organised as follows.

By choosing a positive lower orthant dependent distribution function $C_{\bf Z}$ in the domain of attraction of a max-stable copula $C^*$
(Joe, 1997) we will find a new class of limiting multivariate extreme value (MEV) distribution $H$ for the vector ${\bf
M}_n=\left(M_{n,1},\ldots,M_{n,d}\right)$ of componentwise maxima from ${\bf Y}_1,...,{\bf Y}_n$.

In the third section, we derive the multivariate extremal index of the M5 processes and illustrate the result with some  choices  of $C^*$.

Finally, we compare the bivariate upper tail dependence coefficients of $F_{{\bf Y}}$ with the ones of the limiting MEV distribution $H$.

\section{Domains of max-attraction }

Let $\{{\hat{{\bf Y}}}_n\}_n\geq 1$ be a sequence of independent random vectors associated to $\{{\bf Y}_n\}_n\geq 1$, that is such that $F_{{\hat{{\bf Y}}}_n}=F_{{\bf Y}}$, and  ${\hat {\bf
M}}_n=\left(\hat{M}_{n,1},\ldots,\hat{M}_{n,d}\right)$ be the corresponding vector of pointwise maxima. \\

We shall assume that the copula $C_{\bf Z}$ of the array of innovations in $\{{\bf Y}_n\}_n\geq 1$ is positive lower orthant dependent, that is, it satisfies the inequality 
\begin{eqnarray}\label{plod}
C_{\bf Z}(u_1,...,u_d)\geq \displaystyle\prod_{i=1}^d u_i,\,\,(u_1,...,u_d)\in [0,1]^d.
\end{eqnarray}
Positive lower orthant dependence is implied by several  dependence concepts (Joe, 1997, chapter 2)) and, in particular, is satisfied by associated random variables. The multivariate extreme value distributions are associated and therefore the following results can be applied when the array of innovations ${\bf Z}_{l,n}$ has a MEV distribution.

Following Smith and Weissman (1996), we present in this section the limiting distributions of the normalised vectors  ${\hat {\bf
M}}_n$ and ${\bf M}_n$. \\

\begin{pro}\label{p1.1}
If $C_{\bf Z}$ is positive lower orthant dependent and is in the domain of attraction of a max-stable $C^*$, that is, 
\begin{eqnarray*}\label{p2.1-1}
C_{\bf Z}^n(u_1^{1/n},...,u_d^{1/n}) \conv C^*(u_1,...,u_d),
\end{eqnarray*}
then $C_{\bf Y}$ is in the domain of attraction of
\begin{eqnarray}\label{2a}
{\hat{C}}(u_1,...,u_d)=\prod_{l=1}^{\infty}\prod_{k=-\infty}^{\infty}C^*\left(u_1^{\alpha_{l,k,1}},...,u_d^{\alpha_{l,k,d}}\right)
\end{eqnarray}
\end{pro}

\dem We have 
\vspace{0.3cm}\begin{eqnarray*}
\lim_{n\to \infty} C_{\bf Y}^n(u_1^{1/n},...,u_d^{1/n})
=\lim_{n\to \infty}\prod_{l=1}^{\infty}\prod_{k=-\infty}^{\infty}C_{\bf Z}^n
\left(u_1^{\alpha_{l,k,1}/n},...,u_d^{\alpha_{l,k,d}/n}\right)
\end{eqnarray*}\vspace{0.5cm}
and we can take the product of limits since 
\begin{eqnarray*}
\lim_{n\to \infty}\prod_{l=1}^{\infty}\prod_{k=-\infty}^{\infty}C_{\bf Z}^n
\left(u_1^{\alpha_{l,k,1}/n},...,u_d^{\alpha_{l,k,d}/n}\right)=
\exp\left(\lim_{n\to \infty}\sum_{l=1}^{\infty}\sum_{k=-\infty}^{\infty}\log \left(C_{\bf Z}^n
\left(u_1^{\alpha_{l,k,1}/n},...,u_d^{\alpha_{l,k,d}/n}\right)\right)\right),
\end{eqnarray*}\vspace{0.5cm}
and, by applying the Fr\'echet lower bound in (\ref{plod}), for each $k$ and $l$ it holds
\begin{eqnarray*}
\left|\log C_{\bf Z}^n
\left(u_1^{\alpha_{l,k,1}/n},...,u_d^{\alpha_{l,k,d}/n}\right)\right |\leq 
\left|\log \left(\prod_{j=1}^du_j^{\alpha_{l,k,j}/n}\right)^n\right|\leq
\sum_{j=1}^{d} \alpha_{l,k,j}\left|\log u_j \right|
\end{eqnarray*}

\vspace{0.5cm}
with $\sum_{l=1}^{\infty}\sum_{k=-\infty}^{\infty}\alpha_{l,k,j}\left|\log u_j \right|$ summable for each $j$. It then follows, from the Domminated Convergence Theorem, that
\vspace{0.3cm}\begin{eqnarray*}
\lefteqn{\hspace{0.3cm}\lim_{n\to \infty}\prod_{l=1}^{\infty}\prod_{k=-\infty}^{\infty}C_{\bf Z}^n
\left(u_1^{\alpha_{l,k,1}/n},...,u_d^{\alpha_{l,k,d}/n}\right)}\nonumber\\[0.3cm]
&=&\exp\left(\sum_{l=1}^{\infty}\sum_{k=-\infty}^{\infty}\log \left(\lim_{n\to \infty}C_{\bf Z}^n
\left(u_1^{\alpha_{l,k,1}/n},...,u_d^{\alpha_{l,k,d}/n}\right)\right)\right)={\hat{C}}\left(u_1,...,u_d\right).
\end{eqnarray*}\vspace{0.5cm}
\hspace{15cm}\fdem

 We have then, for 
 $\forall \boldsymbol{\tau}=(\tau_1,\ldots,\tau_d)\in \R^d_+,$ 

 \begin{eqnarray*}
 \lefteqn{\hspace{0.3cm}P\left(\hat{M}_{n,1}\leq \frac{n}{\tau_1},\ldots,\hat{M}_{n,d}\leq \frac{n}{\tau_d}\right)
=C_{\bf Y}^n\left(e ^{-\tau_1/n},\ldots,e ^{- \tau_d/n}\right)}\nonumber\\[0.3cm]
&\conv& \hat{C}\left(e ^{-\tau_1},\ldots,e ^{-\tau_d}\right)= 
 \prod_{l=1}^{\infty}\prod_{k=-\infty}^{\infty}C^*\left(e^{-\tau_1\alpha_{l,k,1}},...,e^{- \tau_d\alpha_{l,k,d}}\right).
\end{eqnarray*}
\vspace{0.5cm}

Let $\hat{H}$ denote the multivariate extreme value distribution with standard Fr\'echet margins and copula $\hat{C}$, i.e. $\hat{H}(x_1,...,x_d)=\displaystyle\prod_{l=1}^{\infty}\prod_{k=-\infty}^{\infty}C^*\left(e^{-x_1^{-1}\alpha_{l,k,1}},...,e^{- x_d^{-1}\alpha_{l,k,d}}\right)$. We will now consider the corresponding limit $H$ for the vector  ${\bf
M}_n=\left(M_{n,1},\ldots,M_{n,d}\right)$ of componentwise maxima from ${\bf Y}_1,...,{\bf Y}_n$.
\vspace{0.5cm}


\begin{pro}\label{p1.2}
If $C_{\bf Z}$ is a positive lower orthant dependent copula in the domain of attraction of a max-stable $C^*$ then
\begin{eqnarray}
\lim_{n\to \infty}P\left(M_{n,1}\leq \frac{n}{\tau_1},\ldots,M_{n,d}\leq \frac{n}{\tau_d}\right)=
\prod_{l=1}^{\infty}C^*\left(e^{-\displaystyle{\max_{-\infty\leq k\leq
+\infty}\alpha_{l,k,1}\tau_1}},...,e^{-\displaystyle{\max_{-\infty\leq k\leq
+\infty}\alpha_{l,k,d}\tau_d}}\right).\label{2}
\end{eqnarray}
\end{pro}

\dem It holds 
\begin{eqnarray*}
\lefteqn{\hspace{0.3cm}\lim_{n\to \infty}P\left(M_{n,1}\leq \frac{n}{\tau_1},\ldots,M_{n,d}\leq \frac{n}{\tau_d}\right)}\nonumber\\[0.3cm]
&=&\lim_{n\to \infty} \prod_{l=1}^{\infty}\prod_{m=-\infty}^{\infty}C_{\bf Z}\left(e^{-\displaystyle{\max_{1-m\leq k\leq
n-m}\alpha_{l,k,1}\tau_1/n}},...,e^{-\displaystyle{\max_{1-m\leq k\leq n-m}\alpha_{l,k,d}\tau_d/n}}\right)\nonumber \\[0.3cm]
&=&\lim_{n\to \infty}\exp\left(\sum_{l=1}^{\infty}\sum_{m=-\infty}^{\infty} \log 
C_{\bf Z}\left(e^{-\displaystyle{\max_{1-m\leq k\leq n-m}\alpha_{l,k,1}\tau_1/n}},...,e^{-\displaystyle{\max_{1-m\leq k\leq n-m}\alpha_{l,k,d}\tau_d/n}}\right)\right).
\end{eqnarray*}

Now we can interchange the first two limits, by using analogous arguments to those used in Proposition 2.1, since

\begin{equation}
\begin{split}
&\sum_{m=-\infty}^{\infty} 
\left|\log C_{\bf Z}\left(e^{-\displaystyle{\max_{1-m\leq k\leq n-m}\alpha_{l,k,1}\tau_1/n}},...,e^{-\displaystyle{\max_{1-m\leq k\leq n-m}\alpha_{l,k,d}\tau_d/n}}\right)\right|\nonumber \\[0.3cm]
&\leq\sum_{m=-\infty}^{\infty} 
\left|\log \left(e^{-\displaystyle{\sum_{j=1}^d\max_{1-m\leq k\leq n-m}\alpha_{l,k,j}\tau_j/n}}\right)\right|\nonumber\\[0.3cm]
&=\sum_{m=-\infty}^{\infty} \frac{1}{n}\displaystyle{\sum_{j=1}^d\max_{1-m\leq k\leq n-m}\alpha_{l,k,j}\tau_j}
\leq\sum_{m=-\infty}^{\infty} \frac{d}{n}\,\displaystyle{\max_{1\leq j\leq d}\,\,\max_{1-m\leq k\leq n-m}\alpha_{l,k,j}\tau_j},
\end{split}
\end{equation}

and this last sum is bounded for all $n$ by $\sum_{k=-\infty}^{\infty} d\displaystyle{\max_{1\leq j\leq d}(\alpha_{l,k,j}\tau_j})$, which is summable in $l$ (Smith and Weissman, 1996).

We have therefore 
\begin{eqnarray*}
\lefteqn{\hspace{0.3cm}\lim_{n\to \infty}P\left(M_{n,1}\leq \frac{n}{\tau_1},\ldots,M_{n,d}\leq \frac{n}{\tau_d}\right)}\nonumber\\[0.3cm]
&=&\exp\left(\sum_{l=1}^{\infty}\lim_{n\to \infty}\sum_{m=-\infty}^{\infty}\frac{1}{n} \log 
C_{\bf Z}^n\left(e^{-\displaystyle{\max_{1-m\leq k\leq n-m}\alpha_{l,k,1}\tau_1/n}},...,e^{-\displaystyle{\max_{1-m\leq k\leq n-m}\alpha_{l,k,d}\tau_d/n}}\right)\right)
\end{eqnarray*}

and we have to prove that 
\begin{eqnarray*}
\lefteqn{\hspace{0.3cm}\lim_{n\to \infty}\sum_{m=-\infty}^{\infty}\frac{1}{n} \log 
C_{\bf Z}^n\left(e^{-\displaystyle{\max_{1-m\leq k\leq n-m}\alpha_{l,k,1}\tau_1/n}},...,e^{-\displaystyle{\max_{1-m\leq k\leq n-m}\alpha_{l,k,d}\tau_d/n}}\right)}\nonumber\\[0.3cm]
&=&\log C^*\left(e^{-\displaystyle{\max_{-\infty\leq k\leq
+\infty}\alpha_{l,k,1}\tau_1}},...,e^{-\displaystyle{\max_{-\infty\leq k\leq
+\infty}\alpha_{l,k,d}\tau_d}}\right)
\end{eqnarray*}

For fixed $l$, let $b_{l,k,j}=a_{l,k,j}\tau _j$, $j=1,...,d$, and $b_{l,k}=\displaystyle{\max_{j=1,...,d}}b_{l,k,j}$ (which are summable in $k$ and $l$). Suppose that $b_{l,k,j}$ is maximized when $k=k^*(l,j)$, $b_{l,k}$ is maximized when $k=k^*(l)$ (not necessarily unique) and assume, without loss of generality, that $k^*(l,1)\leq ...\leq k^*(l,d)$.
Break the above sum into three sums $S_n^{(i)}$, $1=1,2,3$, corresponding to $1-\min\{k^*(l,1), k^*(l)\}\leq m \leq n-\max\{k^*(l,d), k^*(l)\}$, $m< 1-\min\{k^*(l,1), k^*(l)\}$ and $m>n-\max\{k^*(l,d),k^*(l)\}$.

For large $n$ we have 

\begin{eqnarray*}
\lefteqn{\hspace{0.3cm}S_n^{(1)}=\frac{1}{n} \displaystyle\sum_{m=1-\min\{k^*(l,1),k^*(l)\}}^{n-\max\{k^*(l,d), k^*(l)\}} \log 
C_{\bf Z}^n\left(e^{-\displaystyle{\max_{1-m\leq k\leq n-m}\alpha_{l,k,1}\tau_1/n}},...,e^{-\displaystyle{\max_{1-m\leq k\leq n-m}\alpha_{l,k,d}\tau_d/n}}\right)}\nonumber\\[0.3cm]
&=& \frac{n-\max\{k^*(l,d), k^*(l)\}+\min\{k^*(l,1),k^*(l)\}}{n}\log C_{\bf Z}^n\left(e^{-\alpha_{l,k^*(l,1)1}\tau_1/n},...,e^{-\alpha_{l,k^*(l,d)d}\tau_d/n}\right),
\end{eqnarray*}

which converges to 

\begin{eqnarray*}
\log C^*\left(e^{-\alpha_{l,k^*(l,1)1}\tau_1},...,e^{-\alpha_{l,k^*(l,d)d}\tau_d}\right)
=\log C^*\left(e^{-\displaystyle{\max_{-\infty\leq k\leq
+\infty}\alpha_{l,k,1}\tau_1}},...,e^{-\displaystyle{\max_{-\infty\leq k\leq
+\infty}\alpha_{l,k,d}\tau_d}}\right).
\end{eqnarray*}

Otherwise, by applying (\ref{plod}), we can write

\begin{eqnarray*}
\lefteqn{\hspace{0.3cm}\left|S_n^{(2)}\right|\leq \frac{1}{n} \displaystyle\sum_{m<1- k^*(l)}
\displaystyle\sum_{j=1}^d \max_{1-m\leq k\leq n-m}\alpha_{l,k,j}\tau_j}\nonumber\\[0.3cm]
&+&\frac{1}{n} \displaystyle\sum_{1- k^*(l)\leq m < 1-\min\{k^*(l,1), k^*(l)\}}
\displaystyle\sum_{j=1}^d \max_{1-m\leq k\leq n-m}\alpha_{l,k,j}\tau_j\nonumber\\[0.3cm]
&\leq& \frac{d}{n} \displaystyle\sum_{m<1- k^*(l)}
\max_{1-m\leq k\leq n-m}\displaystyle{\max_{j=1,...,d}}\alpha_{l,k,j}\tau_j\nonumber\\[0.3cm]
&+&\frac{d}{n} \displaystyle\sum_{1- k^*(l)\leq m < 1-\min\{k^*(l,1), k^*(l)\}}
\max_{1-m\leq k\leq n-m}\displaystyle{\max_{j=1,...,d}}\alpha_{l,k,j}\tau_j\nonumber\\[0.3cm]
&\leq& \frac{d}{n} \displaystyle\sum_{m<1- k^*(l)}
\max_{1-m\leq k\leq n-m}\displaystyle{\max_{j=1,...,d}}\alpha_{l,k,j}\tau_j+
\frac{d}{n}\left(\min\{k^*(l,1), k^*(l)\}- k^*(l)\right)b_{l,k^*(l)}
\end{eqnarray*}

and 

\begin{eqnarray*}
\lefteqn{\hspace{0.3cm}\left|S_n^{(3)}\right|\leq \frac{1}{n} \displaystyle\sum_{m > n- k^*(l)}
\displaystyle\sum_{j=1}^d \max_{1-m\leq k\leq n-m}\alpha_{l,k,j}\tau_j}\nonumber\\[0.3cm]
&+&\frac{1}{n} \displaystyle\sum_{n- \max \{k^*(l,d),k^*(l)\} < m \leq n-k^*(l)}
\displaystyle\sum_{j=1}^d \max_{1-m\leq k\leq n-m}\alpha_{l,k,j}\tau_j\nonumber\\[0.3cm]
&\leq& \frac{d}{n} \displaystyle\sum_{m > n- k^*(l)}
\max_{1-m\leq k\leq n-m}\displaystyle{\max_{j=1,...,d}}\alpha_{l,k,j}\tau_j\nonumber\\[0.3cm]
&+&\frac{d}{n} \displaystyle\sum_{n- \max \{k^*(l,d),k^*(l)\} < m \leq n-k^*(l)}
\max_{1-m\leq k\leq n-m}\displaystyle{\max_{j=1,...,d}}\alpha_{l,k,j}\tau_j\nonumber\\[0.3cm]
&\leq& \frac{d}{n} \displaystyle\sum_{m > n- k^*(l)}
\max_{1-m\leq k\leq n-m}\displaystyle{\max_{j=1,...,d}}\alpha_{l,k,j}\tau_j+
\frac{d}{n}\left(\max \{k^*(l,d),k^*(l)\}- k^*(l)\right)b_{l,k^*(l)}
\end{eqnarray*}

The second terms in the above upper bounds for $\left|S_n^{(i)}\right|$, $i=1,2$, tend to zero and the same holds for the first terms, by Lemma 3.2 in Smith and Weissman (1996).\\

\hspace{15cm}\fdem
\vspace{0.5cm}

For the particular case of $C^*(u_1,...,u_d)=\displaystyle\min_{1\leq j\leq d} u_j$, $(u_1,\ldots,u_d) \in
[0,1]^d$,  the above result agrees with (3.5) in  Smith and Weissman (1996), that is 
\begin{eqnarray*}
\lim_{n\to \infty}P\left(M_{n,1}\leq \frac{n}{\tau_1},\ldots,M_{n,d}\leq \frac{n}{\tau_d}\right)=
\exp\left\{\displaystyle{-\sum_{l=1}^{\infty}
\max_{-\infty< k<\infty}\max_{1\leq j\leq
d}\alpha_{l,k,j}\tau_j}\right\}.
\end{eqnarray*}

If we choose $C^*(u_1,...,u_d)=\prod_{j=1}^d u_j$, $(u_1,\ldots,u_d) \in
[0,1]^d$, we find again the result in case two of the last example in Martins and Ferreira (2005), 
\begin{eqnarray*}
\lim_{n\to \infty}P\left(M_{n,1}\leq \frac{n}{\tau_1},\ldots,M_{n,d}\leq \frac{n}{\tau_d}\right)=
\exp\left\{\displaystyle{-\sum_{j=1}^{d}\sum_{l=1}^{\infty}
\max_{-\infty< k<\infty}\alpha_{l,k,j}\tau_j}\right\}.
\end{eqnarray*}

In this last work it is considered only $C_{\bf Z}=C^*$ and the above two choices for  $C^*$ are treated directly. 

By calculating the limits in (\ref{2a}) and (\ref{2}) for a given initial copula $C^*$, the above propositions enable us to obtain a larger class of MEV models.

Some choices for $C^*$ are available in the literature. For instance, Cap\'{e}ra\`{a} {it et al.} (2000) gives sufficient conditions for an Archimedean copula $C_{\bf Z}$ to be in the domain of attraction of the Gumbel-Hougaard or logistic copula $C^*(u_1,...,u_d)=\exp \left(-\left(\sum_{j=1}^d\left(-\log u_j\right)^\alpha\right)^{1/\alpha} \right)$. More examples of  $C_{\bf Z}$ in the domain of attraction of a MEV copula $C^*$  can be found in Demarta and McNeil (2005) and H\"{u}sler and Reiss (1989).

\section{The multivariate extremal index of M$5$ process}\label{sdef+prop}

In this section we will extend the results about the multivariate extremal index of the M$4$ process, as a corollary of the propositions 2.1 and 2.2.\\

We first recall the definition of the multivariate extremal index function $\theta(\boldsymbol{\tau})=\theta(\tau_1,\ldots,\tau_d),$
$\boldsymbol{\tau} \in \R^d_+,$ that relates the MEV distribution functions $H$ and
$\hat{H}$ and which was introduced by Nandagopalan (1990).

A $d-$dimensional stationary sequence $\{{\bf Y}_n\}_{n\geq 1}$ is said to have a
mul\-ti\-va\-riate ex\-tremal index $\theta(\boldsymbol{\tau})\in
[0,1]$, $\boldsymbol{\tau} \in \R^d_+,$ if for each $ \boldsymbol{\tau}=(\tau_1,\ldots,\tau_d)$  in 
$\R^d_+,$ there exists $ {\bf u}_n^{(\boldsymbol{\tau})}=$
\linebreak$(u_{n,1}^{(\tau_1)},\ldots,u_{n,d}^{(\tau_d)}),$ $n\geq
1,$ satisfying\\
\begin{eqnarray*} 
    nP(Y_{1j}>u_{n,j}^{(\tau_j)})\conv \tau_j,\,\, j=1,\ldots,d,
\end{eqnarray*}
 \begin{eqnarray*}
 P(\hat{\bf M}_n\leq {\bf u}_n^{(\boldsymbol{\tau})})\conv \hat{\gamma}(\boldsymbol{\tau})\in (0,1]\,\,\,\,and\,\,\,\, P({\bf M}_n\leq {\bf
u}_n^{(\boldsymbol{\tau})})\conv
\hat{\gamma}(\boldsymbol{\tau})^{\theta(\boldsymbol{\tau})}.
\end{eqnarray*}
\vspace{0.3cm}

\begin{pro}\label{p2.1}
If $C_{\bf Z}$  is a positive lower orthant dependent copula in the domain of attraction of a max-stable $C^*$ then the multivariate extremal index  of  the M$5$ process $\{{\bf Y}_{n}\}$  defined in (\ref{1}) is given by 

\begin{eqnarray*}
\theta(\tau_1,\ldots,\tau_d)=\frac{\displaystyle\sum_{l=1}^{\infty}\log C^*\left(e^{-\displaystyle{\max_{-\infty\leq k\leq
+\infty}\alpha_{l,k,1}\tau_1}},...,e^{-\displaystyle{\max_{-\infty\leq k\leq
+\infty}\alpha_{l,k,d}\tau_d}}\right)}
{\displaystyle\sum_{l=1}^{\infty}\displaystyle\sum_{k=-\infty}^{\infty}\log C^*\left(e^{-\alpha_{l,k,1}\tau_1},...,e^{-\alpha_{l,k,d}\tau_d}\right)}
\end{eqnarray*}
and the extremal index of $\{Y_{n,j}\}_{n\geq 1}$ is
$$\theta_{j}=\displaystyle{\sum_{l=1}^{\infty}\max_{-\infty<
k<\infty}\alpha_{l,k,j}},\,\,j=1,\ldots,d.$$
\end{pro}

\dem By applying Proposition 2.1 and Proposition 2.2 with $u_{n,j}^{(\tau_j)}=\frac{\tau_j}{n}$, $j=1,...,d$, we find 

$$\hat{\gamma}(\boldsymbol{\tau})=\prod_{l=1}^{\infty}\prod_{k=-\infty}^{\infty}C^*\left(e^{-\alpha_{l,k,1}\tau_1},...,e^{-\alpha_{l,k,d}\tau_d}\right)$$

and 

$$\gamma(\boldsymbol{\tau})=\hat{\gamma}(\boldsymbol{\tau})^{\theta(\boldsymbol{\tau})}=\prod_{l=1}^{\infty}C^*\left(e^{-\displaystyle{\max_{-\infty\leq k\leq
+\infty}\alpha_{l,k,1}\tau_1}},...,e^{-\displaystyle{\max_{-\infty\leq k\leq
+\infty}\alpha_{l,k,d}\tau_d}}\right),$$

that leads to the results for $\theta(\tau_1,\ldots,\tau_d)=\frac{\log \gamma(\boldsymbol{\tau})}{\log \hat{\gamma}(\boldsymbol{\tau})}$ and $\theta_{j}=\displaystyle\lim_{\tau_i\rightarrow 0^+,\, i\in \{1,...,d\}-\{j\}}\theta(\tau_1,\ldots,\tau_d)$.

\hspace{15cm}\fdem
\vspace{0.3cm}

Therefore the copulas $C^*$, ${\hat{C}}$  and the copula $C$ of the limiting MEV distribution $H$ are related by a multivariate extremal index $\theta(\tau_1,\ldots,\tau_d)$ through

\begin{eqnarray*}\label{Ceta2}
C(u_1,...,u_d)=\left({\hat {C}}(u_1^{1/\theta_1},...,u_d^{1/\theta_d})\right)^{\theta(-\frac{\log u_1}{\theta_1},...,-\frac{\log u_d}{\theta_d})}
=\prod_{l=1}^{\infty}\prod_{k=-\infty}^{\infty}\left(C^*(u_1^{\alpha_{l,k,1}/\theta_1},...,u_d^{\alpha_{l,k,d}/\theta_d})\right)^{\theta(-\frac{\log u_1}{\theta_1},...,-\frac{\log u_d}{\theta_d})}.
\end{eqnarray*}

The proposition 3.1 leads to the results of Smith and Weissman (1996) and Martins and Ferreira (2005) when  $C^*$ is the copula of  the minimum and product, respectively. 

Otherwise, if we take for instance the Logistic copula $C^*$ we find 
$$\theta(\tau_1,\ldots,\tau_d)=\frac
{\displaystyle\sum_{l=1}^{\infty}
\left(\sum_{j=1}^{d}
\left(\displaystyle{\max_{-\infty\leq k\leq+\infty}\alpha_{l,k,j}\tau_j}\right)^\alpha\right)^{1/\alpha}}
{\displaystyle\sum_{l=1}^{\infty}\sum_{k=\infty}^{\infty}\left(\sum_{j=1}^{d}
\left(\alpha_{l,k,j}\tau_j\right)^{\alpha} \right)^{1/\alpha}}.$$

\vspace{0.5cm}

  


\section{The tail dependence of M$5$ process}\label{sdef+prop}

For a random vector ${\bf X}=(X_1,\ldots,X_d)$ with continuous margins $F_1,...,F_d$ and copula $C$, let the bivariate (upper) tail dependence coefficients parameters be defined by 

\begin{eqnarray*}
\lambda_{j,j'}^{(C)} =\lim_{u\uparrow 1} P(F_j(X_j)>u|F_{j'}(X_{j'})>u),\,\,1\leq j<j'\leq d.
\end{eqnarray*}

The tail dependence coefficient characterizes the dependence in the tail of a random pair
$(X_j,X_{j'})$, i.e., $\lambda_{j,j'}^{(C)}>0$ corresponds to tail dependence and
$\lambda_{j,j'}^{(C)}=0$ means tail independence, and can be
defined via the copula of the random vector which refers to their
dependence structure  independently of
their marginal distributions.

It holds 
\begin{eqnarray*}
 \lambda_{j,j'}^{(C)} =2-\lim_{u\uparrow 1} \frac{\log C_{j,j'}(u,u)}{\log u},
\end{eqnarray*}

where $C_{j,j'}$ is the copula of the sub-vector $(X_j,X_{j'})$. 

In this section we will relate $\lambda_{j,j'}^{(C)}$ with $\lambda_{j,j'}^{({\hat{C}})}$, where $C$ is copula of $H$ and ${\hat{C}}$ is the copula of ${\hat{H}}$.
We first remark that, for each $(u_j,u_{j'}) \in [0,1]^2$, we have 
\begin{eqnarray}\label{3}
C_{j,j'}(u_j,u_{j'})=\left({\hat {C}}_{j,j'}(u_j^{1/\theta_j},u_{j'}^{1/\theta_{j'}})\right)^{\theta(-\frac{\log u_j}{\theta_j},-\frac{\log u_{j'}}{\theta_{j'}})},
\end{eqnarray}
where $\theta(\tau_j,\tau_{j'})$ is the bivariate extremal index of $\{(Y_{n,j},Y_{n,j'})\}_{n\geq 1}$.
This relation enables us to compare the tail dependence parameters $\lambda_{j,j'}^{(C)}$ with $\lambda_{j,j'}^{({\hat{C}})}$ through the function $\theta(\tau_j,\tau_{j'})$.
 
 \begin{pro}\label{p3.1}
If $C$ and ${\hat{C}}$ satisfy (\ref{3}) then
\begin{itemize}
\item[(a)]  
$\lambda_{j,j'}^{(C)}=2+\theta(\frac{1}{\theta_j},\frac{1}{\theta_{j'}}) \log {\hat {C}}_{j,j'}(e^{-1/\theta_j},e^{-1/\theta_{j'}}),$

\item[(a)]  $\lambda_{j,j'}^{(C)}=\lambda_{j,j'}^{({\hat{C}})}+\log \frac{ {\hat {C}}_{j,j'}(e^{-\theta(\frac{1}{\theta_j},\frac{1}{\theta_{j'}})/\theta_j},e^{-\theta(\frac{1}{\theta_j},\frac{1}{\theta_{j'}})/\theta_{j'}})}{ {\hat {C}}_{j,j'}(e^{-1},e^{-1})}.$
\end{itemize}
\end{pro}
 
 \dem
  From the spectral measure representation (Resnick, 1987) of the copula ${\hat{C}}$, with measure ${\hat{W}}$, we get 
   \begin{eqnarray}\label{4}
\lambda_{j,j'}^{({\hat{C}})} =2-\lim_{u\uparrow 1} \frac{\log {\hat{C}}_{j,j'}(u,u)}{\log u}=2-\int_{\mathcal{S}_d} 
\max\{w_j,w_{j'}\}\, d{\hat{W}}(w_1,...,w_d)=2+\log {\hat{C}}_{j,j'}(e^{-1},e^{-1}).
\end{eqnarray}
 By using (\ref{3}) and the homogeneity of order 0 of the multivariate extremal index, it follows that 
 \begin{eqnarray*}
 \lefteqn{\hspace{0.3cm}\lambda_{j,j'}^{(C)} =2-\lim_{u\uparrow 1} \frac{\log C_{j,j'}(u,u)}{\log u}=
 2-\lim_{u\uparrow 1}\theta(\frac{-\log u}{\theta_j},\frac{-\log u}{\theta_{j'}}) 
  \frac{\log {\hat{C}}_{j,j'}(u_j^{1/\theta_j},u_{j'}^{1/\theta_{j'}})}{\log u}}\nonumber\\[0.3cm]
&=&2-\theta(\frac{1}{\theta_j},\frac{1}{\theta_{j'}})\lim_{u\uparrow 1}
\frac{\int_{\mathcal{S}_d} \max\{\frac{-\log u w_j}{\theta _j},\frac{-\log u w_{j'}}{\theta _{j'}}\}
\, d{\hat{W}}(w_1,...,w_d)}
{-\log u} \nonumber\\[0.3cm]
&=&2-\theta(\frac{1}{\theta_j},\frac{1}{\theta_{j'}})
\int_{\mathcal{S}_d} \max\{\frac{w_j}{\theta _j},\frac{ w_{j'}}{\theta _{j'}}\}
\, d{\hat{W}}(w_1,...,w_d) \nonumber\\[0.3cm]
&=&2+\theta(\frac{1}{\theta_j},\frac{1}{\theta_{j'}})\log {\hat {C}}_{j,j'}(e^{-1/\theta_j},e^{-1/\theta_{j'}}),
\end{eqnarray*}
which has (\ref{4}) as a particular case.
To obtain the second statement we combine the first with (\ref{4}) and use the max-stability of ${\hat{C}}_{j,j'}$, as follows:
\begin{eqnarray*}
 \lefteqn{\hspace{0.3cm}\lambda_{j,j'}^{(C)} =2+\log {\hat{C}}_{j,j'}(e^{-1},e^{-1})-\log {\hat{C}}_{j,j'}(e^{-1},e^{-1})
+ \theta(\frac{1}{\theta_j},\frac{1}{\theta_{j'}})\log {\hat {C}}_{j,j'}(e^{-1/\theta_j},e^{-1/\theta_{j'}})}\nonumber\\[0.3cm]
 &=&\lambda_{j,j'}^{({\hat{C}})}+\log \frac
 {\left({\hat {C}}_{j,j'}(e^{-1/\theta_j},e^{-1/\theta_{j'}})\right)^{\theta(\frac{1}{\theta_j},\frac{1}{\theta_{j'}})}}
 {{\hat{C}}_{j,j'}(e^{-1},e^{-1})}\nonumber\\[0.3cm]
 &=&\lambda_{j,j'}^{({\hat{C}})}+\log \frac
 {{\hat {C}}_{j,j'}(e^{-\theta(\frac{1}{\theta_j},\frac{1}{\theta_{j'}})/\theta_j},e^{-\theta(\frac{1}{\theta_j},\frac{1}{\theta_{j'}})/\theta_{j'}})}
 {{\hat{C}}_{j,j'}(e^{-1},e^{-1})}.
\end{eqnarray*}

\hspace{15cm}\fdem\\

We will now apply these relations to the particular case of M$5$ processes, enhancing the effect of the multivariate extremal index in the tail dependence.

 \begin{pro}\label{p3.2}
Let $\{{\bf Y}_n\}_n\geq 1$ be a M$5$ process defined as in (\ref{1}) and such that  $C_{\bf Z}$ a positive lower orthant dependent copula in the domain of attraction of $C^*$. Then, for any $1\leq j <j'\leq d$ it holds 
\begin{itemize}
\item[(a)]  $\lambda_{j,j'}^{({\hat {C}})}=2+ \displaystyle\sum_{l=1}^{\infty}\displaystyle\sum_{k=-\infty}^{\infty}
\log C_{j,j'}^{*}\left(e^{-\alpha_{l,k,j}},e^{-\alpha_{l,k,{j'}}}\right).$
\end{itemize}
\begin{itemize}
\item[(b)]  $\lambda_{j,j'}^{(C)}=2+\theta(\frac{1}{\theta_j},\frac{1}{\theta_{j'}}) \displaystyle\sum_{l=1}^{\infty}\displaystyle\sum_{k=-\infty}^{\infty}
\log C_{j,j'}^{*}\left(e^{-\alpha_{l,k,j}/\theta_j},e^{-\alpha_{l,k,{j'}}/\theta_{j'}}\right).$
\end{itemize}
\end{pro}

\dem By using (\ref{4}) and the Proposition 2.1, we obtain (a); by using the Proposition (4.1)-(a) and Proposition 2.1, we get (b).
\hspace{11cm}\fdem\\
\vspace{0.5cm}

If $C^*(u_1,...,u_d)=\displaystyle\min_{1\leq j\leq d}u_j$, $(u_1,\ldots,u_d) \in
[0,1]^d$, we obtain 
\begin{eqnarray}\label{5}
\lambda_{j,j'}^{({\hat {C}})}=2- \displaystyle\sum_{l=1}^{\infty}\displaystyle\sum_{k=-\infty}^{\infty}
\max\{\alpha_{l,k,j},\alpha_{l,k,{j'}}\}
\end{eqnarray}

 in general greater than zero, which agrees with the result (2.10) presented in  Heffernan {\it et al.} (2007).
We can also say from (b) that, for this copula,
$$\lambda_{j,j'}^{(C)}=2- \displaystyle\sum_{l=1}^{\infty}
\max \left\{\displaystyle{\max_{-\infty\leq k\leq
+\infty}}\alpha_{l,k,j}/\theta_j, \displaystyle{\max_{-\infty\leq k\leq
+\infty}}\alpha_{l,k,{j'}}/\theta_{j'}\right\}$$

and it holds $\lambda_{j,j'}^{(C)}>\lambda_{j,j'}^{({\hat {C}})}$ if and only if 
$$ \displaystyle\sum_{l=1}^{\infty}\displaystyle\sum_{k=-\infty}^{\infty}
\max\{\alpha_{l,k,j}, \alpha_{l,k,{j'}}\}- \displaystyle\sum_{l=1}^{\infty}\max\left\{
\displaystyle{\max_{-\infty\leq k\leq
+\infty}}\alpha_{l,k,j}/\theta_j, \displaystyle{\max_{-\infty\leq k\leq
+\infty}}\alpha_{l,k,{j'}}/\theta_{j'}\right\}>0,$$
where $\theta_{j}=\displaystyle{\sum_{l=1}^{\infty}\max_{-\infty<
k<\infty}\alpha_{l,k,j}},\,\,j=1,\ldots,d.$

We can easily construct examples of M4 processes for which $\lambda_{j,j'}^{(C)}>\lambda_{j,j'}^{({\hat {C}})}$ or
$\lambda_{j,j'}^{(C)}<\lambda_{j,j'}^{({\hat {C}})}$.

The results for this class of moving multivariate maxima show us that if we obtain or estimate a tail dependence parameter of the common distribution of the variables in a stationary sequence we don't have necessarily the corresponding parameter in the limiting MEV model $C$. The multivariate extremal index of the stationary sequence can increase or decrease the tail dependence and the extremal coefficients of the limiting MEV model ${\hat {C}}$ arising from the i.i.d.  sequence.

\vspace{0.5cm}

Since $\lambda_{j,j'}^{({\hat {C}})}=\lambda_{j,j'}^{(C_{{\bold Y}})}$, the result in (\ref{5}) says that in the M4 processes the  variables $Y_{n,j}$, $j=1,...,d$, are in general asymptotically dependent. For the M5 processes we can choose $C^*$ in order to produce variables asymptotically independent. Take for instance  $C^*$ such that $\lambda_{j,j'}^{(C^*)}=0$ and $\alpha_{l,k,j}=\alpha_{l,k,j'}$. We then find $\lambda_{j,j'}^{({\hat {C}})}=\lambda_{j,j'}^{(C^*)}=0$.

\vspace{0.5cm}

The Proposition 4.2. also points out that even for a choice of a copula $C^*$ with symmetric tail dependencies, the values of the signatures $\alpha_{l,k,j}$ can lead to asymmetric tail dependencies in both  copulas $C$ and ${\hat {C}}$, that is, to different values of $\lambda_{j,j'}$ for different choices of $(j,j')$.
\vspace{0.5cm}

The results in the above proposition are translations of the classical result $\lambda=2-\epsilon$ for $C_{j,j'}$ and ${\hat {C}}_{j,j'}$ 
(Nelsen, 2006)
In fact, if we define the extremal coefficient of the MEV copula $C$ as the constant $\epsilon_C$ such that $C(u,...,u)=u^{\epsilon_C}$ for all $u\in [0,1]$, then from 
$$ C(u_1,...,u_d)=\exp\bigg(-\theta(\frac{\log u_1}{\theta_1},...,\frac{\log u_d}{\theta_d})
\int_{\mathcal{S}_d} \max_{1\leq j \leq d} \frac{-\log u_j w_j}{\theta_j}\, d{\hat{W}}(w_1,...,w_d)\bigg),$$
we find
$$\epsilon_C=-\theta(\frac{1}{\theta_1},...,\frac{1}{\theta_d})\log {\hat {C}}\left(e^{-1/\theta_1},...,e^{-1/\theta_d}\right)$$
and, in particular,  
$$\epsilon_ {\hat {C}}=-\log {\hat {C}}\left(e^{-1},...,e^{-1}\right).$$
The relation between ${\hat {C}}_{j,j'}$ and $C^*_{j,j'}$ enables now to obtain the Proposition 4.2. from $\lambda_{j,j'}^{({\hat {C}})}=2-\epsilon_ {{\hat {C}}_{j,j'}}$ and  $\lambda_{j,j'}^{(C)}=2-\epsilon_{C_{j,j'}}$.
\vspace{0.5cm}

Other features of the M5 processes such as a directory of tail dependence coefficients for different $C^*$ and signatures $\alpha_{l,k,j}$, illustrating the range of dependence structures, the model selection and estimation are  key directions in a future research.
\vspace{0.5cm}



\biblio{

\bibitem{caperaa+} Cap\'{e}ra\`{a}, P., Foug\`{e}res, A.L. and Genest, C. (2000)  Bivariate distributions with given extreme value attractor. J. Multivariate Anal., 72, 30-49.

\bibitem{davis+resnick} Davis, R.A. and Resnick, S.I. (1989) Basic properties and prediction of Max-ARMA processes. Adv. Appl. Prob., 21, 781-803.

\bibitem{deheuvels} Deheuvels, P. (1983) Point processes and multivariate extreme values. J. Multivariate Anal., 13, 257-252.

\bibitem{Demarta+McNeil} Demarta, S. and McNeil, A. (2005) The t-copula and related copulas. International Statistical Review 73, 111-129.

\bibitem{hall+} Hall, P., Peng, L. and Yao, Q. (2002) Moving-maximum models for extrema of time series. J. Stat. Plan. Infer., 103, 51-63.

\bibitem{hef+} Heffernan, J. E., Tawn, J. A., Zhang, Z. (2007). Asymptotically
(in)dependent multivariate maxima of moving maxima processes,
Extremes, 10, 57-82.

\bibitem{Husler+Reiss} H\"{u}sler, J. and Reiss, R. (1989) Maxima of normal random vectors: between independence and complete dependence. Statist. Probab. Lett., 76, 283-286.

\bibitem{joe} Joe, H. (1997). Multivariate Models and Dependence Concepts. Chapman
\& Hall, London.

\bibitem{martins+ferreira} Martins, A.P. and Ferreira, H. (2005) The extremal index and the dependence structure of a multivariate extreme value distribution. Test, Vol. 14, 2, 433-448.

\bibitem{nand} Nandagopalan, S. (1990). Multivariate extremes and estimation of the
extremal index. Ph.D. Dissertation, Department of Statistics,
University of North Carolina, USA.

\bibitem{nelsen_cop} Nelsen, R.B. (2006). An Introduction to Copulas. Second Edition.
Springer, New York.

\bibitem{res} Resnick, S. (1987). \emph{Extreme Values, Regular Variation and Point Processes.} Springer-Verlag, New York.

\bibitem{smith+weissman} Smith, R.L. and Weissman, I. (1996). Characterization and estimation of the
multivariate extremal index. Technical Report, Univ. North
Carolina.

\bibitem{zhang+smith} Zhang, Z. and Smith, R.L. (2004) The behavior of multivariate maxima of moving maxima processes. J. Appl. Probab., 41, 1113-1123.

}

\end{document}